\documentclass[12pt]{amsart}
\usepackage{amsmath,amssymb,txfonts}
\usepackage{amssymb}
\usepackage{amsmath}
\usepackage{mathrsfs}
\usepackage{amsmath,amssymb}

\newtheorem{theorem}{Theorem}[section]
\newtheorem{lemma}[theorem]{Lemma}

\theoremstyle{definition}
\newtheorem{definition}[theorem]{Definition}

\theoremstyle{remark}
\newtheorem{remark}[theorem]{Remark}

\def\R{\mbox{\rm{I}}\!\mbox{\rm{R}}}

\def\<{\leq}             \def\>{\geq}

\numberwithin{equation}{section}

\usepackage{color}

\begin{document}

\title[Navier-Stokes equations]
{A remark on ill-posedness}


\author[Haibo Yang]{Haibo Yang}
\address{Faculty of Mathematics and Statistics, Hubei University,
Wuhan, 430062, China.}
\email{yanghb97@qq.com}

\author[Qixiang Yang]{Qixiang Yang}
\address{
School of Mathematics and Statistics, Wuhan University, Wuhan, 430072, China.
}
\email{qxyang@whu.edu.cn}

\author[Huoxiong Wu]{Huoxiong Wu}
\address{School of Mathematical sciences, Xiamen University, Xiamen Fujian, 361005, China.
}
\email{huoxwu@xmu.edu.cn}


\thanks{This work was supported by the National Natural Science Foundation of China (No. 11571261, 11771358, 11871101).}


\date{}

\dedicatory{}


\begin{abstract}
Norm inflation implies certain discontinuous dependence of the solution on the initial value.
The well-posedness of the mild solution means the existence and uniqueness of the fixed points of the corresponding integral equation.
For ${\rm BMO}^{-1}$, Auscher-Dubois-Tchamitchian \cite{ADT} proved that Koch-Tataru's solution is stable.
In this paper, we construct a non-Gauss flow function to show that, for classic Navier-Stokes equations,
wellposedness and norm inflation may have no conflict and stability may have meaning different to $L^{\infty}(({\rm BMO}^{-1})^{n})$.
\end{abstract}

\maketitle


 \vspace{0.1in}

\section{Introduction and main results} 
\label{intro}
In this paper, we give a remark on the meaning of convergence of mild solution and the illposedness of the following incompressible Navier-Stokes equations:
\begin{equation}\label{eqn:ns}
\left\{\begin{array}{ll}  \partial_{ t} u
+(-\Delta)^{\beta} u + u \cdot \nabla u -\nabla p=0,
& \mbox{ in }  [0,T)\times \mathbb{R}^{n}; \\
\nabla \cdot u=0,
& \mbox{ in } [0,T)\times \mathbb{R}^{n}; \\
u|_{t=0}= u_0, & \mbox{ in } \mathbb{R}^{n};
\end{array}
\right.
\end{equation}
where $u(t,x)$ and $p(t,x)$ denote the velocity vector field and the pressure of fluid at the point $(t,x)\in [0,T)\times \mathbb{R}^{n}$ respectively. While $u_0$ is a given initial velocity vector field.
The wellposedness for different initial data spaces have been studied heavily. See Cannone \cite{C1}, Iwabuchi-Nakamura \cite{IN}, Koch-Tataru \cite{KT}, Li-Xiao-Yang \cite{LXY}, Yang-Yang \cite{YY}.
The solutions of the above Cauchy problem can be obtained via the integral equation:
\begin{equation}\label{eqn:mildsolution}
u(t,x)= e^{-t(-\Delta)^{\beta} } u_0(x) - B(u,u)(t,x),
\end{equation}
where
\begin{equation}\label{eqn:HW}
\begin{cases}
B(u,u)(t,x)\equiv\int^{t}_{0} e^{-(t-s)(- \Delta)^{\beta}}
\mathbb{P}\nabla (u\otimes u) ds,\\
\mathbb{P}\nabla (u\otimes u)\equiv \sum\limits_{l}
\partial x_{l} (u_{l}u) -
\sum\limits_{l} \sum\limits_{l'} (-\Delta)^{-1} \partial x_{l}\partial x_{l'} \nabla (u_{l} u_{l'}).
\end{cases}
\end{equation}
The equation (\ref{eqn:mildsolution}) can be solved by a fixed-point method whenever the convergence
is suitably defined in certain function spaces.
For $u_0$ belongs to some initial space $X^n= (X(\mathbb{R}^{n}))^{n}$, denote
\begin{equation}\label{eqn:it}
\begin{cases}
u^{(0)}(t,x) = e^{-t(-\Delta)^{\beta}} u_0,\\
u^{(\tau +1)} (t,x) = u^{(0)}(t,x) - B(u^{(\tau)}, u^{(\tau)})(t,x), \forall \tau=0,1,2,\cdots,
\end{cases}
\end{equation}
where $e^{t\Delta}u_0$ belongs to some space $Y^{n}= (Y((0,T)\times \mathbb{R}^n))^{n}$.

The above iteration process convergence for $\|u_0\|_{X^n}$ small enough.
Such solutions of (\ref{eqn:mildsolution}) are called mild solutions of
(\ref{eqn:ns}). The notion of such a mild solution was pioneered by
Kato-Fujita \cite{KF} in 1960s. During the latest decades, many
important results about mild solutions to (\ref{eqn:ns}) have been established.
Given $t\in (0,T]$ and $u(t,x)$ belongs to the Banach space $X^{n}$.
We know $u(t,x)$ belongs to function space $L^{\infty}((0,T],X^{n})$ means $$(\|u(t,x)\|_{X^{n}})_{L^{\infty}((0,T])}<\infty.$$
For initial data $u_0\in X^{n}$, most often, its solution $u(t,x)$ belong to some solution spaces $Y(X^{n})$ which is a subspace of $L^{\infty}((0,T],X^{n})$.
See, for example, Cannone \cite{C1, C2},
Germin-Pavlovic-Staffilani \cite{GPS}, Giga-Miyakawa \cite{GM}, Kato \cite{Kat},  Lemari\'e \cite{Lem, Lem1},
Wu \cite{W1,W2,W3,W4} and some author's collaboration work \cite{LXY, LiY, LY,YL}.

When $\beta=1$, the above equations \eqref{eqn:ns} are the classic Navier-Stokes equations.
When a solution $u(t,x)$ is not bounded in $L^{\infty}((0,T],X^{n})$,
Bourgain-Pavlovi\'c \cite{BP} and Yoneda \cite{Yo} call such solution has norm inflation phenomenon and
for corresponding initial value spaces $X^{n}$,
the equations are ill-posed in the  sense of $L^{\infty}((0,T],X^{n})$ norm.
For the end point Triebel-Lizorkin spaces $(\dot{F}^{-1,q}_{\infty})^{n}(2<q\leq \infty)$,
Bourgain-Pavlovi\'c \cite{BP} and Yoneda \cite{Yo} have shown the norm inflation in the end Triebel-Lizorkin spaces $L^{\infty}((0,T], (\dot{F}^{-1,q}_{\infty})^{n})$.
Wang \cite{Wa} has shown norm inflation in the end point Besov spaces $(\dot{B}^{-1, q}_{\infty})^{n}$.
But all these do not mean that the mild solution cannot have other strong stability.

Koch-Tataru introduced ${\rm BMO}^{-1}$ in Theorem 1 of \cite{KT}.
We know $ \dot{F}^{-1,2}_{\infty}={\rm BMO}^{-1}$. See \cite{LY}, \cite{Me} and \cite{YSY}.
${\rm BMO}^{-1}$ has a very special structure, and it has inspired a lot of interest in studying it.
Let $l(Q)$ denote the side length of cube $Q$.
Auscher-Dubois-Tchamitchian \cite{ADT} introduced $u(t,x)\in C_{T}$.
\begin{definition}\label{def:BMOsolution}
$u(t,x)\in C_{T}$ if and only if the following two conditions are satisfied:
\begin{equation}\label{eq:infty}
N_{\infty,T}(u)
{\Tiny\begin{array}{c}
def\\ =
\end{array}}
\sup\limits_{0<t\leq T} t^{\frac{1}{2}} \| u(t,x)\|_{\infty}<\infty.
\end{equation}
\begin{equation}\label{eq:BMO-1}
N_{c,T}(u){\Tiny\begin{array}{c}
def\\ =
\end{array}}
\sup\limits_{\mbox{ cube } Q: l^2(Q) \leq T}
\{|Q|^{-1} \int^{l^2(Q)}_{0}\int_{Q} |u(t,x)|^{2} dtdx\}^{\frac{1}{2}}<\infty.
\end{equation}
If $u(t,x)$ satisfies further
$$\lim\limits_{T'\rightarrow 0} \|u\|_{C_{T'}}=0,$$
then $u(t,x)\in C_{0,T}$. Denote $C_0= C_{0,\infty}$.
\end{definition}
Auscher-Dubois-Tchamitchian \cite{ADT} proved that
Koch-Tataru's solution is stable in the sense of $C_{0}$ for initial data in ${\rm VMO}^{-1}$
where $C_{0}$ has nothing to do with $L^{\infty}({\rm BMO}^{-1})^{n}$.
See Remark \ref{re:3.2}.
Miura \cite{Mi} has made a remark on uniqueness of mild solutions to the
Navier-Stokes equations for initial data in $({\rm BMO}^{-1}\bigcap L^{\infty})^{n}$.
Germin-Pavlovic-Staffilani \cite{GPS} considered the regularity of solution with initial data in $({\rm BMO}^{-1})^{n}$.
Therefore, when we discuss the relationship between ${\rm BMO}^{-1}$ and the well-posedness,
it is   necessary to clarify the relationship between $L^{\infty}((0,T], ({\rm BMO}^{-1})^{n})$ and Koch-Tataru's solution space in \cite{KT}.

To establish the wellposedness for initial data in  $({\rm BMO}^{-1})^{n}$,
Koch-Tataru \cite{KT} introduced the following space $Y((0,T],({\rm BMO}^{-1})^{n})$:
\begin{definition}\label{def:BMOsolution}
$u(t,x)\in Y((0,T],({\rm BMO}^{-1})^{n})$ if and only if $u(t,x)\in C_{T}$ and the following condition is satisfied:
\begin{equation}\label{eq:zero}
\nabla u(t,x)=0 {\mbox { in } } [0,T]\times \mathbb{R}^{n},
\end{equation}
\end{definition}
For initial data in $({\rm BMO}^{-1})^{n}$,
Koch-Tataru \cite{KT} have established wellposedness in solution space $Y((0,T],({\rm BMO}^{-1})^{n})$.
We find Koch-Tataru's space allows norm inflation in $L^{\infty}((0,T],({\rm BMO}^{-1})^{n})$.
In fact, we prove that Koch-Tataru's space $Y((0,T],({\rm BMO}^{-1})^{n})$ is not a subspace of
$L^{\infty}((0,T],\linebreak ({\rm BMO}^{-1})^{n})$.
That is, the well-posedness and norm inflation have no conflict,\linebreak
Auscher-Dubois-Tchamitchian's stability result shows
stability can mean different things. Precisely, we obtain the following theorem.
\begin{theorem}\label{th:BMOsolution}
Koch-Tataru's space $Y((0,1],({\rm BMO}^{-1})^{n})$ is not a subspace of
$L^{\infty} ((0,1], ({\rm BMO}^{-1})^{n}) $.
That is, there exists
$u(t,x)\in Y((0,1],({\rm BMO}^{-1})^{n})$ but
$\|u(t,x)\|_{L^{\infty} ((0,1],({\rm BMO}^{-1})^{n})}=\infty$.
\end{theorem}
The rest of this paper is organized as follows: In section 2, we will present some preliminaries about Meyer wavelets,
then we present wavelet characterization for end point Triebel-Lizorkin spaces and Koch-Tataru's solution space.
In section 3, we use Meyer wavelets to construct some functions in Koch-Tataru's space and prove Theorems \ref{th:BMOsolution}.

\section{Wavelets and function spaces}\label{sec2}

In this section, we recall first some auxiliary knowledge on wavelets.
We indicate that we will use tensorial product real valued
orthogonal Meyer wavelets. We refer the reader to \cite{Me, Woj, Yang1} for further information.
Let $\Psi^{0}$ be an even function in $ C^{\infty}_{0}
([-\frac{4\pi}{3}, \frac{4\pi}{3}])$ with
\begin{equation}
\left\{ \begin{aligned}
&0\leq\Psi^{0}(\xi)\leq 1; \nonumber\\
&\Psi^{0}(\xi)=1\text{ for }|\xi|\leq \frac{2\pi}{3}.\nonumber
\end{aligned} \right.
\end{equation}
  Write
  $$\Omega(\xi)= \sqrt{(\Psi^{0}(\frac{\xi}{2}))^{2}-(\Psi^{0}(\xi))^{2}}.$$ Then $\Omega(\xi)$ is an even
function in $ C^{\infty}_{0}([-\frac{8\pi}{3}, \frac{8\pi}{3}])$.
Clearly,

\begin{equation}
\left\{ \begin{aligned}
&\Omega(\xi)=0\text{ for }|\xi|\leq \frac{2\pi}{3};\nonumber\\
&\Omega^{2}(\xi)+\Omega^{2}(2\xi)=1=\Omega^{2}(\xi)+\Omega^{2}(2\pi-\xi)\text{
for }\xi\in [\frac{2\pi}{3},\frac{4\pi}{3}].\nonumber
\end{aligned} \right.
\end{equation}

 Let $\Psi^{1}(\xi)=
\Omega(\xi) e^{-\frac{i\xi}{2}}$. For any $\epsilon=
(\epsilon_{1},\cdots, \epsilon_{n}) \in \{0,1\}^{n}$, define
$\Phi^{\epsilon}(x)$ by $\hat{\Phi}^{\epsilon}(\xi)=
\prod\limits^{n}_{i=1} \Psi^{\epsilon_{i}}(\xi_{i})$. For $j\in
\mathbb{Z}$ and $k\in\mathbb{Z}^{n}$, let $\Phi^{\epsilon}_{j,k}(x)=
2^{\frac{nj}{2}} \Phi^{\epsilon} (2^{j}x-k)$. $\forall \epsilon\in \{0,1\}^{n}, j\in \mathbb{Z}, k\in \mathbb{Z}^{n}$ and distribution $f(x)$, denote
$f^{\epsilon}_{j,k}=\langle f, \Phi^{\epsilon}_{j,k}\rangle$.
Furthermore, we put
\begin{equation}\nonumber
\Lambda_{n} =\{(\epsilon,j,k), \epsilon\in \{0,1\}^{n}\backslash\{0\}, j\in\mathbb{Z},
k\in \mathbb{Z}^{n}\}.
\end{equation}

Sobolev space $\dot{H}^{\frac{n}{2}-1}= \dot{F}^{\frac{n}{2}-1,2}_{2}$,
Lebesgue space $L^{n}= \dot{F}^{0,2}_{n}$,
Besov spaces $\dot{B}^{\frac{n}{p}-1,p}_{p}=\dot{F}^{\frac{n}{p}-1,p}_{p}$
and ${\rm BMO}^{-1}= \dot{F}^{-1,2}_{\infty}$ are all Triebel-Lizorkin spaces.
For an overview of function spaces, we refer to Li-Xiao-Yang \cite{LXY},
Lin-Yang \cite{LY}, Yang \cite{Yang1} and Yuan-Sickel-Yang \cite{YSY}.
Denote $\mathfrak{D}=\{Q_{j,k}=2^{-j}k+2^{-j}[0,1]^{n}, \forall j\in \mathbb{Z}, k\in \mathbb{Z}^{n}\}$.
We recall then the wavelet characterization of end-point Triebel-Lizorkin spaces $\dot{F}^{\gamma,q}_{\infty}(\mathbb{R}^{n})$ (see \cite{LXY,LY,YSY}).
\begin{lemma}\label{lem:c}
Given $1\leq q \leq \infty$ and $\gamma\in \mathbb{R}$.
$f(x)= \sum\limits_{\epsilon,j,k} a^{\epsilon}_{j,k}
\Phi^{\epsilon}_{j,k}(x)\in
\dot{F}^{\gamma,q}_{\infty}(\mathbb{R}^{n})\Leftrightarrow$
\begin{equation}\label{eq:bmc}
\begin{array}{rl}
&\sup\limits_{Q\in \mathfrak{D} } \Big\{|Q|^{-1}
\sum\limits_{(\epsilon,j,k):Q_{j,k}\subset Q}
2^{jq(\gamma_{1}+\frac{n}{2}-\frac{n}{q})}
|a^{\epsilon}_{j,k}|^{q}\Big\}^{\frac{1}{q}} <+
\infty.\end{array}
\end{equation}

\end{lemma}

At the end of this section, we present one lemma on Koch-Tataru's solution space.
Let $u(t,x)= (u_1(t,x), u_2(t,x),\dots, u_n(t,x))^{t}$. For $i=1,2,\cdots, n$,
denote $$u_{i}(t,x) = \sum\limits_{(\epsilon,j,k)\in \Lambda_n} a^{i,\epsilon}_{j,k}(t) \Phi^{\epsilon}_{j,k}(x).$$
The following lemma  is a direct corollary of the wavelet characterization in Lemma \ref{lem:c}.
\begin{lemma}\label{le:4.1}
(i) $u(t,x)$ satisfies \eqref{eq:BMO-1} if and only if
\begin{equation}\label{B.BMO}
\sup\limits_{i, j_0\geq 0, k_0\in \mathbb{Z}^{n}} 2^{nj_0} \int^{2^{-2j_0}}_{0} \sum\limits_{Q_{j,k}\subset Q_{j_0, k_0}}
|a^{i,\epsilon}_{j,k}(t)|^{2} dt <\infty.
\end{equation}

(ii) $u(t,x)\in L^{\infty} ((0,1], ({\rm BMO}^{-1})^{n})$ if and only if
\begin{equation}\label{B.lim}
\sup\limits_{i, 0< t\leq 1, k_0\in \mathbb{Z}^{n}} 2^{nj_0}  \sum\limits_{Q_{j,k}\subset Q_{j_0, k_0}}
2^{-2j} |a^{i,\epsilon}_{j,k}(t)|^{2}  <\infty.
\end{equation}
\end{lemma}

\section{Proof of Theorem \ref{th:BMOsolution}}
We give first two remarks on the mild solution.
The wellposedness of mild solution is independent of whether the initial value space is separable or not.
\begin{remark} \label{re:3.1}
(i) The unique existence of the mild solution $u(t,x)$ of integral equation \eqref{eqn:mildsolution} means
$u_{0}$ belongs to some initial space $X^{n}$ and there exists some solution space $Y((0,T], X^{n})$ such that
\begin{itemize}
\item [(1)] $e^{-t (-\Delta)^{\beta}} u_0 \in Y((0,T], X^{n})$.
\item [(2)] The iteration process \eqref{eqn:it} converges to $u(t,x)\in Y((0,T], X^{n})$.
\end{itemize}
Auscher-Dubois-Tchamitchian \cite{ADT} believes that
the stability of solution $u(t,x)$ refers to convergence to the boundary value $u_{0}$
in the sense of norm matching the solution space.

(ii) Whether or not the solution space $Y((0,T], X^{n})$ is a subspace of $L^{\infty}( X^{n})$
is independent of whether or not $X^{n}$ is separable. For fractional Navier-Stokes equations where $\frac{1}{2}<\beta<1$,
Yu-Zhai \cite{YZ} considered the non-separable Bloch spaces $\dot{B}^{1-2\beta, \infty}_{\infty}(\mathbb{R}^{n})$.
Hence such fractional equations have both wellposedness and continuous dependence.
\end{remark}

For classic Navier-Stokes equations,
according to the above Theorem \ref{th:BMOsolution},
norm inflation and wellposedness are not incompatible.
Let $X^{n}$ be the initial data space.
The Wellposedness does not require solution to be contained in $L^{\infty} ((0,T], X^{n})$.
\begin{remark}\label{re:3.2}
(i) Bourgain-Pavlovi\'c  and Yoneda's illposedness
results mean only the norm inflation in $L^{\infty} ((0,T], X^{n})$,
don't mean that it's unstable in any other sense.
Wellposedness means that the iterative process \eqref{eqn:it}
can converge to a particular and unique function in certain particular space.
This supports Chemin and Gallagher's point in \cite{CG}:
the well-posedness need not the boundedness in $L^{\infty} ((0,T], X^{n})$.

(ii) Germin-Pavlovic-Staffilani's regularity result in \cite{GPS}
simply indicates that the solution is sufficiently smooth, independent of stability.
When Miura consider the stability of ${\rm BMO}^{-1}$ in \cite{Mi},
he added the condition that the initial value is in the intersection with $L^{\infty}$.
Miura considered $({\rm BMO}^{-1}\bigcap L^{\infty})^{n}$
and the corresponding stability shows only the stability under the restriction of $L^{\infty}$.

(iii) Auscher-Dubois-Tchamitchian \cite{ADT} proved that
Koch-Tataru's solution is stable in the sense of $C_{0}$ for initial data in ${\rm VMO}^{-1}$.
This is the strongly convergent boundary behavior of the solution in the sense of solution space norm.
Auscher-Dubois-Tchamitchian's space $C_{0,T}$ is defined for non Gauss flow.
Our result shows that a function $u(t,x)$ belongs to $C_{0}$ has nothing to do with
that $u(t,x)$ belongs to $L^{\infty}([0,\infty), ({\rm BMO}^{-1})^{n})$.
Auscher-Dubois-Tchamitchian's stability provides another kind of stability different to $L^{\infty}( X^{n})$
norm sense.
\end{remark}

Now we come to prove Theorem \ref{th:BMOsolution}.
\begin{proof}
Denote $\Lambda= \{(e,j,k)\in \Lambda_n,  e=(1,\cdots,1), j\geq 0, k\in \mathbb{Z}^{n}\}$.
Take $0<a<\frac{1}{2}$ and $\frac{n}{2} +2a-1<b <\frac{n}{2}$ and
take $u_{1}(t,x) = \sum\limits_{(e,j,k)\in \Lambda} a^{e}_{j,k}(t) \Phi^{e}_{j,k}(x)$ where $a^{e}_{j,k}(t)$ satisfies
$$a^{e}_{j,k}(t)= \left\{
\begin{array}{lll}
t^{-a} 2^{-bj}, & 1\leq j \leq -\frac{1}{2} \log_2 t, &k\in \mathbb{Z}^{n};\\
0, & j > -\frac{1}{2} \log_2 t {\mbox { or }} t\geq 1; &k\in \mathbb{Z}^{n}.
\end{array}
\right.$$
We know, if $0<b<\frac{n}{2}$ and $b\geq \frac{n}{2}+  2a-1$, then $u_{1}(t,x)$ satisfies the following equation:
\begin{equation*}
t^{\frac{1}{2}} \|u_1(t,x)\|_{\infty}<\infty.\end{equation*}

The number of $k$ satisfying $Q_{j,k}\subset Q_{j_0, k_0}$ is $2^{n(j-j_0)}$.
Hence for $j_0\geq 0$,
$$\begin{array}{rl}
&2^{nj_0} \int^{2^{-2j_0}}_{0} \sum\limits_{Q_{j,k}\subset Q_{j_0, k_0}} |a^{e}_{j,k}(t)|^{2} dt\\
\leq & 2^{nj_0} \int^{2^{-2j_0}}_{0} t^{-2a} \sum\limits_{j\geq j_0, 1\leq j \leq -\frac{1}{2} \log_2 t} 2^{-2bj} 2^{n(j-j_0)} dt\\
= &  \int^{2^{-2j_0}}_{0} t^{-2a} \sum\limits_{\max (j_0, 1)\leq j \leq -\frac{1}{2} \log_2 t} 2^{(n-2b)j}  dt\\
\leq C &  \int^{2^{-2j_0}}_{0} t^{b-\frac{n}{2}-2a}   dt.
\end{array}$$
If $b> \frac{n}{2}+  2a-1$, then $u_{1}(t,x)$ satisfies equation \eqref{B.BMO}.

If $ a<\frac{1}{2}$ and $\frac{n}{2} +2a-1<b <\frac{n}{2}$, then
$$\begin{array}{rl}
c_{t,j_0}=& 2^{nj_0}  \sum\limits_{Q_{j,k}\subset Q_{j_0, k_0}}
2^{-2j} |a^{e}_{j,k}(t)|^{2} \\
\leq & 2^{nj_0}  t^{-2a} \sum\limits_{j\geq j_0, 1\leq j \leq -\frac{1}{2} \log_2 t} 2^{-2j-2bj} 2^{n(j-j_0)} \\
= &   t^{-2a} \sum\limits_{\max (j_0, 1)\leq j \leq -\frac{1}{2} \log_2 t} 2^{(n-2b-2)j}.
\end{array}$$
Take $j_0=0, b>\frac{n}{2}-1$ and $a>0$, then $c_{t,0}> c t^{-2a}$, hence $\sup\limits_{0<t<1} c_{t,0}=+\infty$.
Hence, $u_{1}(t,x)$ does satisfy equation \eqref{B.lim}.

Take $u_{2}(t,x)= -\frac{1}{\partial_2} \partial_1 u_1(t,x)$. It is easy to see that $u_{2}(t,x)$ satisfies the same properties as
$u_{1}(t,x)$ does.

For $i=3,\cdots, n$, take $u_{i}(t,x)=0$. By construction, we know that $u(t,x)= \linebreak (u_1(t,x), u_2(t,x), \cdots, u_n(t,x))^{t}$ satisfies
$$\nabla u(t,x)=0 {\mbox { in }} [0,1]\times \mathbb{R}^{n}.$$
That is, $u(t,x)$ satisfies all the conditions in Theorem \ref{th:BMOsolution}.
\end{proof}

{\bf Acknowledgements.}
The authors would like to thank Professor Baoxiang Wang for his invaluable discussions and suggestions.
Further, the corresponding author would like to thank Professor Zhifei Zhang and Zihua Guo for their invaluable discussions and suggestions.


\end{document}